\documentclass[a4paper,11pt,twoside]{article}
\RequirePackage[a4paper,head=1cm]{geometry}
\usepackage[utf8x]{inputenc}
\usepackage{textcomp}
\usepackage{amsmath,amsfonts,verbatim,afterpage,euscript,mathrsfs,amssymb,empheq}
\usepackage{amsthm}
\usepackage{dsfont}
\usepackage{pst-fill,pst-grad}
\usepackage{textcomp}
\usepackage{multicol}
\usepackage{amsmath}
\usepackage{amssymb}
\usepackage{hyperref}
\usepackage{dsfont}
\usepackage{colortbl}

\newtheorem{Theoreme}{Theorem}

\newtheorem{Corollaire}{Corollary}
\newtheorem{Proposition}{Proposition}

\newcommand{\mysection}{\setcounter{equation}{0} \section}

\title{\bf  Two pointwise estimates for singular integral operators over $\mathfrak{p}$-adic groups} 

\author{Diego Chamorro\footnote{Laboratoire de Math\'ematiques et Mod\'elisation d'Evry (LaMME) - UMR 8071. Universit\'e d'Evry Val d'Essonne, 23 Boulevard de France, 91037 Evry Cedex, France. email: \textit{diego.chamorro@univ-evry.fr}}.} 

\begin{document} 
\maketitle 
\abstract{We establish some new pointwise boundedness properties for rough singular integral operators $T_\Omega$ defined over $\mathfrak{p}$-adic groups $\mathbb{Q}^n_\mathfrak{p}$. The properties of these operator rely essentially on the behavior of their rough kernels $\Omega$ and we will present two different cases that fall outside the setting of the classical Calder\'on-Zygmund setting.\\

\noindent\textbf{Keywords:}  singular integral operator; $\mathfrak{p}$-adic groups.\\
\textbf{MSC 2020}: 22E35; 43A15.
}

\section{Introduction}
In this article we study some pointwise estimates for rough singular integral operators $T_\Omega$ defined by the expression 
\begin{equation}\label{Def_Operator_intro}
T_\Omega(f)(x)=\underset{\epsilon\to 0}{\lim}\int_{\{\epsilon<\|x-y\|_{\mathfrak{p}}\}}\Omega(x,y)f(y)dy,
\end{equation}
where $f:\mathbb{Q}^n_{\mathfrak{p}}\longrightarrow \mathbb{R}$ is a locally measurable function defined over the $n$-dimensional local field $\mathbb{Q}^n_{\mathfrak{p}}$ and where $\Omega:\mathbb{Q}^n_{\mathfrak{p}}\times \mathbb{Q}^n_{\mathfrak{p}}\longrightarrow \mathbb{R}$ is the kernel associated to the operator $T_\Omega$.\\

The main purpose of this article is to prove, under some mild assumption on the kernel $\Omega$, some generic pointwise estimates of the form 
\begin{equation}\label{Pointwise_Estimate1}
|T_\Omega(f)(x)|\leq C \mathscr{M}_\alpha(f)^{\theta}(x)\times \|f\|_{\mathcal{M}^{p,n+\sigma}}^{1-\theta},
\end{equation}
where $0<\theta<1$ is a power index, where for a parameter $0<\alpha<n$, $\mathscr{M}_\alpha(f)$ is the fractional maximal function of $f$ and it is defined by the expression
\begin{equation}\label{Def_FuncMax}
\mathscr{M}_\alpha(f)(x)=\underset{r>0}{\sup} \frac{r^\alpha}{|B(x,r)|}\int_{B(x,r)}|f(y)|dy,
\end{equation}
and where the quantity $\|\cdot\|_{\mathcal{M}^{p,n+\sigma}}$ (with $1\leq p<+\infty$ and $0<\sigma<p$) that characterizes the membership to the Morrey-Campanato space $\mathcal{M}^{p,n+\sigma}(\mathbb{Q}^n_{\mathfrak{p}})$ is defined by the semi-norm
\begin{equation}\label{Def_MorreyCampanato}
\|f\|_{\mathcal{M}^{p,n+\sigma}}=\underset{r>0}{\sup}\underset{x\in \mathbb{Q}^n_{\mathfrak{p}}}{\sup} \left(\frac{1}{r^{n+\sigma}}\int_{B(x,r)}|f(y)-f_B|^pdy\right)^{\frac{1}{p}}.
\end{equation}

\noindent In the setting of local fields, operators of the form (\ref{Def_Operator_intro}) with $\Omega(x,y)=\Omega(x-y)$ (\emph{i.e.} convolution type operators) where studied in \cite{Phillips} under a regularity hypothesis over the kernel $\Omega$ (see also \cite{Phillips0} for the particular case of the Hilbert transform). The main results stated in these references are the boundedness of such operators $T_\Omega$ (under some constraints on the kernels $\Omega$) in the usual Lebesgue spaces:
$$\|T_{\Omega}(f)\|_{L^p}\leq C\|f\|_{L^p},$$
whenever we have $1<p<+\infty$.\\

Let us mention now that more recently, in \cite{Ash}, the boundedness of some similar operators of the type $T_\Omega$ presented above was studied in the framework of Triebel-Lizorkin and Besov spaces over $\mathfrak{p}$-adic groups while Morrey-Herz spaces were considered in \cite{Shi}. See also \cite{Volosivets0}, \cite{Volosivets} and \cite{Wu} for some results about the boundedness of some classical harmonic analysis operators in the setting of $\mathbb{Q}^n_{\mathfrak{p}}$. Remark that in all these references, the operators considered are slightly different from the ones studied here, indeed, we consider in this article a generalization of the classical Calder\'on-Zygmund operators since the properties of the kernels $\Omega$ examined in the lines below do not fall into the usual setting (see \cite{Grafakos} for classical results on these operators in the framework of the Euclidean space).\\

In order to state our results, we need to be more precise with the general framework used here. To begin with, let us fix one prime number $\mathfrak{p}\geq 2$ and let us consider the $n$-dimensional $\mathfrak{p}$-adic group $\mathbb{Q}^n_{\mathfrak{p}}=\mathbb{Q}_{\mathfrak{p}}\times \cdots \times \mathbb{Q}_{\mathfrak{p}}$ for $n\geq 2$. The elements of $\mathbb{Q}^n_{\mathfrak{p}}$ will be denoted by $x=(x_1, \cdots, x_n)$. Note that the group $\mathbb{Q}^n_{\mathfrak{p}}$ can be endowed with an ultra-metric distance $d(x,y)=\|x-y\|_{\mathfrak{p}}$ and with a Haar measure $dx$, in particular we will denote by $|E|$ the measure of the set $E\subset \mathbb{Q}^n_{\mathfrak{p}}$. See Section \ref{IntroP} below for more details about the framework considered here.\\

\noindent Our first theorem is the following one:
\begin{Theoreme}\label{Theorem1}
Consider a locally integrable function $f:\mathbb{Q}^n_{\mathfrak{p}}\longrightarrow \mathbb{R}$ such that $f$ belongs to the Morrey-Campanato space $\mathcal{M}^{1,n+\sigma}(\mathbb{Q}^n_{\mathfrak{p}})$ defined in the expression (\ref{Def_MorreyCampanato}) above with $0<\sigma<1$.\\

\noindent Consider a singular integral operator $T_\Omega$ of the form given in the formula (\ref{Def_Operator_intro}) and assume that its kernel $\Omega(\cdot, \cdot):\mathbb{Q}_{\mathfrak{p}}^n\times \mathbb{Q}_{\mathfrak{p}}^n\longrightarrow \mathbb{R}$ satisfies the following properties:
\begin{itemize}
\item[1)] for all $x,y\in  \mathbb{Q}_{\mathfrak{p}}^n$ such that $x\neq y$, we have the estimate
\begin{equation}\label{SizeKernel}
|\Omega(x,y)|\leq \frac{C_\Omega}{\|x-y\|_{\mathfrak{p}}^n},
\end{equation}
\item[2)] for all $x\in \mathbb{Q}_{\mathfrak{p}}^n$ and for any real numbers $0<a<b<+\infty$, we have the cancellation property
\begin{equation}\label{PseudoRadial}
\int_{\{a<\|x-y\|_{\mathfrak{p}}<b\}}\Omega(x,y)dy=0.
\end{equation}
\end{itemize}
Under these hypotheses, we have the following pointwise control:
\begin{equation}\label{PointwiseControl}
|T_\Omega(x)|\leq C\mathscr{M}_{\alpha}(f)^\theta(x)\|f\|_{\mathcal{M}^{1,n+\sigma}}^{1-\theta},
\end{equation}
where $\mathscr{M}_{\alpha}$ is the fractional maximal function defined in the expression (\ref{Def_FuncMax}) above with $0<\alpha<n$ and the parameter $0<\theta<1$ is fixed by the condition $\theta=\frac{\sigma}{\alpha+\sigma}$.\\
\end{Theoreme}
Let us first remark here that the condition (\ref{SizeKernel}) allows a possible blow up of the kernel $\Omega$ when $x=y$ while the condition (\ref{PseudoRadial}) is crucial in order to give a sense to the integral limit (\ref{Def_Operator_intro}) that defines the operator $T_\Omega$. These two conditions are quite natural in the general setting of singular integral operators (see the book \cite[Section 8.1.1]{Grafakos}). Note however, that besides the conditions (\ref{SizeKernel}) and (\ref{PseudoRadial}), no regularity conditions are asked to the kernel $\Omega$ and this more general framework seems to be new when working with $\mathfrak{p}$-adic spaces and it is an improvement of the results stated in \cite{Phillips}. Observe that similar estimates were studied in \cite{CMM} in the Euclidean case.\\

Remark now that the pointwise control over the operator $T_\Omega$ given in (\ref{PointwiseControl}) is actually given by a mix of two informations: in the one hand we have a fractional maximal function $\mathscr{M}_\alpha$ which will control the large scales (as we shall see in the proof of the Theorem \ref{Theorem1}), while on the other hand, we have a Morrey-Campanato space $\mathcal{M}^{1,n+\sigma}$ that will control the small scales and these two informations are needed to deduce the wished pointwise estimate.\\

Note that the fractional maximal function $\mathscr{M}_\alpha$ conveys some integrability properties in the sense that it usually has a similar smoothing properties as the Riesz potentials (see \emph{e.g.} the formula (\ref{MaxFracRiesz}) below) and remark that the Morrey-Campanato spaces $\mathcal{M}^{p, n+\sigma}$ describe a certain form of regularity (H\"older-type regularity to be more precise, see \cite{Campanato}). In this sense, the pointwise estimate (\ref{PointwiseControl}) can be viewed as a compromise between regularizing effects (via the fractional maximal function) and some regularity constraints (via the Morrey-Campanato spaces).\\

To continue, let us note that we do not know if it is possible to set $\alpha=0$ in the parameter $0<\alpha<n$ that governs the fractional maximal function $\mathscr{M}_\alpha$ (defined in the expression (\ref{Def_FuncMax}) above). Remark that if we set $\alpha=0$ we recover the usual Hardy-Littlewood function, but it seems that the properties of this classical maximal function are not enough to perform our computations.\\

The previous result can be generalized if we have more information over the kernel $\Omega$. Indeed, if we assume that we have $\Omega(x,y)=\Omega(x-y)=\frac{\omega (x-y)}{\|x-y\|_{\mathfrak{p}^n}}$ where the function $\omega$ satisfies a cancellation property as well as some integrability conditions, then Theorem \ref{Theorem1} can be improved as follows:
\begin{Theoreme}\label{Theorem2}
Consider a kernel $\Omega$ of the form
\begin{eqnarray}
\Omega:\mathbb{Q}_{\mathfrak{p}}^n&\longrightarrow &\mathbb{R}\notag\\
z&\longmapsto& \Omega(z)=\frac{\omega(z)}{\|z\|_{\mathfrak{p}}^n},\label{Kernel2}
\end{eqnarray}
where the function $\omega$ satisfies the following cancellation property
\begin{equation}\label{PseudoRadial1}
\int_{\{a<\|z\|_{\mathfrak{p}}<b\}}\frac{\omega(z)}{\|z\|_{\mathfrak{p}}^n}dz=0,
\end{equation}
for all $0<a<b<+\infty$. Assume that $\omega \in L^\rho(\mathbb{Q}^{n}_{\mathfrak{p}})$ for some $1<\rho\leq+\infty$ and define a parameter $\rho'$ such that $1\leq \rho'<+\infty$ by the relationship $1=\frac{1}{\rho}+\frac{1}{\rho'}$. Consider now a singular integral operator $T_\Omega$ given by the expression 
\begin{equation}\label{Def_Operateur2}
T_\Omega(f)(x)=\underset{\epsilon\to 0}{\lim}\int_{\{\epsilon<\|x-y\|_{\mathfrak{p}}\}}\frac{\omega(x-y)}{\|x-y\|_{\mathfrak{p}}^n}f(y)dy.
\end{equation}
If $f:\mathbb{Q}^n_{\mathfrak{p}}\longrightarrow \mathbb{R}$ is a locally integrable function such that $f$ belongs to the Morrey-Campanato space $\mathcal{M}^{\rho',n+\sigma}(\mathbb{Q}^n_{\mathfrak{p}})$, with $1\leq \rho'\leq n+\sigma<+\infty$ and with $0<\sigma< \rho'$.\\

\noindent If we assume moreover the constraint $\omega\in L^{\rho}(\mathbb{Q}^n_\mathfrak{p})$ with $\frac{n+\sigma}{n+\sigma-1}\leq \rho\leq +\infty$, then we have the pointwise bound 
\begin{equation}\label{PointwiseControl1}
|T_{\Omega}(f)(x)|\leq C[\mathscr{M}_{\alpha}(|f|^{\rho'})(x)]^{\frac{\theta}{\rho'}}\|f\|_{\mathcal{M}^{\rho', n+\sigma}}^{1-\theta},
\end{equation}
where $\mathscr{M}_\alpha$ with $0<\alpha<n$ is the fractional maximal function given in (\ref{Def_FuncMax}) and where $0<\theta<1$ is such that $\theta=\frac{\frac{(n+\sigma-n\rho')}{\rho'}}{\frac{(n+\sigma-n\rho')}{\rho'}+\frac{1}{\rho}+\frac{\alpha}{\rho'}}$.\\
\end{Theoreme}
\noindent Let us explain first in which sense this last result is a generalization of the one presented in the Theorem \ref{Theorem1}. Indeed, we remark that if $\omega\in L^\infty(\mathbb{Q}^n_\mathfrak{p})$, (\emph{i.e.} in the particular case when $\rho=+\infty$ and $\rho'=1$), then from the expression (\ref{Kernel2}) we easily deduce the control 
$$|\Omega(x-y)|=\left|\frac{\omega(x-y)}{\|x-y\|_{\mathfrak{p}}^n}\right|\leq \frac{\|\omega\|_{L^\infty}}{\|x-y\|_{\mathfrak{p}}^n}=\frac{C_\omega}{\|x-y\|_{\mathfrak{p}}^n},$$
which is the condition (\ref{SizeKernel}) asked in the Theorem \ref{Theorem1}. Moreover, when $\rho=+\infty$ and $\rho'=1$, we have $\theta=\frac{\frac{(n+\sigma-n\rho')}{\rho'}}{\frac{(n+\sigma-n\rho')}{\rho'}+\frac{1}{\rho}+\frac{\alpha}{\rho'}}=\frac{\sigma}{\alpha+\sigma}$, which is exactly the formula for the index $\theta$ of the Theorem \ref{Theorem1}. Finally note that when $\rho'=1$ we have $\mathcal{M}^{\rho', n+\sigma}=\mathcal{M}^{1, n+\sigma}$ and as we can see in this very particular case we recover \emph{exactly} the statement of the Theorem \ref{Theorem1}. We note now that the case when $\omega\in L^\rho(\mathbb{Q}^n_\mathfrak{p})$ with $\rho<+\infty$ can not be treated directly with the Theorem \ref{Theorem1} and these remarks explain in which sense Theorem \ref{Theorem2} is a generalization of our first result.\\

Let us note now that there is a constraint in the integrability condition of the function $\omega$ as we are asking here $\omega \in L^{\rho}$ with $\rho\leq +\infty$ and $\frac{n+\sigma}{n+\sigma-1}\leq \rho$. The upper bound gives, passing to the index $\rho'$, the condition $1\leq \rho'$ with is natural as we are going to compute some (local) $L^{\rho'}$ norms while the lower bound gives the conditions $\rho'\leq n+\sigma$ and this condition is necessary to define some meaningful Morrey-Campanato spaces.\\

Remark also that if the parameters $\rho'$ and $\sigma$ are big (recall that we have the conditions $\rho'<n+\sigma$ and $\sigma <\rho'<+\infty$), then the index $\rho$ can be close to $1$ (since we have $\frac{n+\sigma}{n+\sigma-1}\leq \rho$), but we can not have $\rho=1$. An interesting open problem is to study if these estimates remain if we replace the limit space $L^1$ -which can not be reached following our method- by the smaller space $L\mbox{log}L$.\\

Finally, as a last observation, it is important to remark that, as before, no regularity is asked to the kernels $\Omega$ (nor to the functions $\omega$). \\

Let us turn now to some consequences of the previous results. In that sense we have the following corollary: 
\begin{Corollaire}\label{Coro1}
Consider the singular integral operator $T_\Omega$ defined by the expression (\ref{Def_Operateur2}) above and assume that the general hypotheses of the Theorem \ref{Theorem2} are in force (in particular consider the same constraints for the parameters $\alpha, \rho, \rho', \sigma$ and $\theta$).\\

\noindent If $1<q<+\infty$ is a real parameter such that $q>\frac{\rho'}{\theta}$ and such that $\frac{n\rho'q\theta}{n\rho'+\alpha q \theta}>1$, if $f$ belongs to the Lebesgue space $L^p(\mathbb{Q}^n_\mathfrak{p})$ with $p=\frac{n\rho'q\theta}{n\rho'+\alpha q \theta}$ and $f$ belongs to the Morrey-Campanato space $\mathcal{M}^{\rho', n+\sigma}(\mathbb{Q}^n_\mathfrak{p})$, then we have the following inequality:
$$\|T_{\Omega}(f)\|_{L^q}\leq C\|f\|_{L^{p}}^{\theta}\|f\|_{\mathcal{M}^{\rho', n+\sigma}}^{1-\theta}.$$
\end{Corollaire}
\noindent Let us remark that this estimate, as well as the pointwise inequalities (\ref{PointwiseControl}) and (\ref{PointwiseControl1}), seem to be new (to the best of our knowledge) in the setting of $\mathfrak{p}$-adic groups when a rough singular integral operator of the form $T_\Omega$ is considered.\\

The plan of this article is the following. In Section \ref{IntroP} we recall some important facts about $\mathfrak{p}$-adic groups, both in terms of their inner topological structure and the functional spaces defined upon them. In Section \ref{ProofTheorem1} we prove Theorem \ref{Theorem1} and Section \ref{ProofTheorem2} is devoted to the proof of the Theorem \ref{Theorem2}. Finally, in Section \ref{ProofCoro1} we present the proof of the Corollary \ref{Coro1}.
\mysection{$\mathfrak{p}$-adic groups, generalities and functional spaces}\label{IntroP}
In this section we recall some well known facts that will be needed here. For more details, see the books \cite{Amice, Koblitz, VVZ} and the references therein. In what follows, $\mathfrak{p}$ will always denote a prime number.\\ 

First, let us write $a|b$ when $a$ divide $b$ or, equivalently, when $b$ is a multiple of $a$. For $0\neq x\in \mathbb{Z}$, we define the $\mathfrak{p}$-adic valuation of $x$ by $\gamma(x)=\max\{r: \mathfrak{p}^{r}|x\}\geq 0$ and, for any rational number $x=\frac{a}{b}\in \mathbb{Q}$, we write $\gamma(x)=\gamma(a)-\gamma(b)$. Furthermore if $x=0$, we agree to write $\gamma(0)=+\infty$. Consider now $x \in \mathbb{Q}$ and let $\mathfrak{p}$ be any prime number, with the $\mathfrak{p}$-adic valuation of $x$ we can construct a norm by writing
\begin{equation}\label{Normeadique}
|x|_{\mathfrak{p}}=\left\lbrace \begin{array}{ll}
\mathfrak{p}^{-\gamma} &\mbox{if}\quad x\neq 0,\\[3mm]
\mathfrak{p}^{-\infty}=0 &\mbox{if}\quad x= 0.
\end{array}\right.
\end{equation}
This expression satisfy the following properties
\begin{enumerate}
\item[1)]$|x|_{\mathfrak{p}}\geq 0$, and $|x|_{\mathfrak{p}}=0 \iff x=0$;
\item[2)] $|xy|_{\mathfrak{p}}=|x|_{\mathfrak{p}}|y|_{\mathfrak{p}}$;
\item[3)] $|x+y|_{\mathfrak{p}}\leq \max\{|x|_{\mathfrak{p}},|y|_{\mathfrak{p}}\}$, with equality when $|x|_{\mathfrak{p}}\neq |y|_{\mathfrak{p}}$.
\end{enumerate}
When a norm satisfy $3)$ it is called a non-Archimedean norm and an interesting fact is that over $\mathbb{Q}$ \textit{all} the possible norms are equivalent to $|\cdot|_{\mathfrak{p}}$ for some $\mathfrak{p}$: this is the so-called Ostrowski theorem, see \cite{Amice} or \cite[Theorem 1, p. 3]{Koblitz} for a proof of this fact.\\

\noindent We now define the field of $\mathfrak{p}$-adic numbers $\mathbb{Q}_{\mathfrak{p}}$ as the completion of $\mathbb{Q}$ when using the norm $|\cdot|_{\mathfrak{p}}$.\\ 

\noindent We present in the following lines the algebraic structure of the set $\mathbb{Q}_{\mathfrak{p}}$. Every $\mathfrak{p}$-adic number $x\neq 0$ can be represented in a unique manner by the formula
\begin{equation}\label{Hensel}
x=\mathfrak{p}^{\gamma}(x_{0}+x_{1}\mathfrak{p}+x_{2}\mathfrak{p}^{2}+...),
\end{equation}
where $\gamma=\gamma(x)$ is the $\mathfrak{p}$-adic valuation of $x$ and $x_{j}$ are integers such that $x_{0}>0$ and $0\leq x_{j}\leq \mathfrak{p}-1$ for $j=1,2,...$. Remark that this canonical representation implies the identity $|x|_{\mathfrak{p}}=\mathfrak{p}^{-\gamma}$.\\

Let $x,y \in \mathbb{Q}_{\mathfrak{p}}$, using the formula (\ref{Hensel}) we define the sum of $x$ and $y$ by $x+y=\mathfrak{p}^{\gamma (x+y)}(c_{0}+c_{1}\mathfrak{p}+c_{2}\mathfrak{p}^{2}+...)$  with $0\leq c_{j}\leq \mathfrak{p}-1$ and $c_{0}>0$, where $\gamma(x+y)$ and $c_{j}$ are the unique solution of the equation
$$\mathfrak{p}^{\gamma(x)}(x_{0}+x_{1}\mathfrak{p}+x_{2}\mathfrak{p}^{2}+...)+ \mathfrak{p}^{\gamma(y)}(y_{0}+y_{1}\mathfrak{p}+y_{2}\mathfrak{p}^{2}+...)=\mathfrak{p}^{\gamma(x+y)}(c_{0}+c_{1}\mathfrak{p}+c_{2}\mathfrak{p}^{2}+...).$$
Furthermore, for $a,x \in \mathbb{Q}_{\mathfrak{p}}$, the equation $a+x=0$ has a unique solution in $\mathbb{Q}_{\mathfrak{p}}$ given by $x=-a$. In the same way, the equation $ax=1$ has a unique solution in $\mathbb{Q}_{\mathfrak{p}}$: $x=1/a$.\\

\noindent We take now a closer look at the topological structure of $\mathbb{Q}_{\mathfrak{p}}$. With the norm $|\cdot|_\mathfrak{p}$ we construct a distance over $\mathbb{Q}_\mathfrak{p}$ by writing
\begin{equation}\label{dist}
d(x,y)=|x-y|_{\mathfrak{p}} 
\end{equation} 
and we define the balls $B_{\gamma}(x)=\left\{y\in \mathbb{Q}_{\mathfrak{p}}: \; d(x,y)\leq \mathfrak{p}^{\gamma}\right\}$  with $\gamma \in \mathbb{Z}$. Remark that, from the properties of the $\mathfrak{p}$-adic valuation, this distance has the \textit{ultra-metric} property (\textit{i.e.} $d(x,y)\leq \max\{d(x,z),d(z,y)\}\leq |x|_{\mathfrak{p}}+|y|_{\mathfrak{p}}$).\\

\noindent We gather with the next proposition some important facts concerning the balls in $\mathbb{Q}_{p}$.
\begin{Proposition}
Let $\gamma$ be an integer, then we have
\begin{enumerate}
\item[1)] the ball $B_{\gamma}(x)$ is a open and a closed set for the distance (\ref{dist}),
\item[2)] every point of $B_{\gamma}(x)$ is its center,
\item[3)] $\mathbb{Q}_{\mathfrak{p}}$ endowed with this distance is a complete Hausdorff metric space,
\item[4)] $\mathbb{Q}_{\mathfrak{p}}$ is a locally compact set,
\item[5)] the $\mathfrak{p}$-adic group $\mathbb{Q}_{\mathfrak{p}}$ is a totally discontinuous space.
\end{enumerate}
\end{Proposition}
\noindent See a proof of this result in \cite[Chapter I]{Koblitz}, see also the books \cite{Amice} and \cite{VVZ} for more details concerning the topological structure of the $p$-adic groups.\\

\noindent Now, for an integer $n\geq 2$ and for a prime number $\mathfrak{p}$, we define the space 
$$\mathbb{Q}_\mathfrak{p}^n=\underbrace{\mathbb{Q}_\mathfrak{p}\times\cdots \times \mathbb{Q}_\mathfrak{p}}_{n},$$ 
and over this group we define the norm $\|\cdot\|_{\mathfrak{p}}$ by the expression
$$\|\cdot\|_{\mathfrak{p}}=\underset{1\leq j\leq n}{\max}|x_j|_{\mathfrak{p}}.$$
Since $\mathbb{Q}_{\mathfrak{p}}^n$ is a locally compact commutative group, there exists a Haar measure $dx$ which is translation invariant \emph{i.e.} $d(x+a)=dx$ and we will normalize the measure $dx$ by setting
$$\int_{\{\|x\|_{\mathfrak{p}}\leq1\}}dx=1.$$
This measure is then unique and we will note $|E|$ the measure for any subset $E$ of $\mathbb{Q}_{\mathfrak{p}}^n$. Remark in particular that we have, for any $k\in \mathbb{Z}$:
$$|B(x,\mathfrak{p}^k)|=\int_{\{\|x\|_{\mathfrak{p}}\leq \mathfrak{p}^k\}}dx=\mathfrak{p}^{nk},$$
and 
$$\int_{\{\mathfrak{p}^{k-1}\leq \|x\|_{\mathfrak{p}}\leq \mathfrak{p}^k\}}dx=\mathfrak{p}^{nk}-\mathfrak{p}^{n(k-1)}.$$
For $1\leq p< +\infty$ the Lebesgue spaces $L^{p}(\mathbb{Q}_{\mathfrak{p}}^n)$ are thus defined in a natural way by asking the condition 
$$\|f\|_{L^p}=\left(\displaystyle{\int_{\mathbb{Q}_{\mathfrak{p}}^n}}|f(x)|^{p}dx\right)^{\frac{1}{p}}<+\infty,$$ 
with the usual modifications when $p=+\infty$. If we have $1=\frac{1}{p}+\frac{1}{p'}$ with $1\leq p, p'\leq +\infty$, then for any measurable subset $K$ of $\mathbb{Q}^n_{\mathfrak{p}}$ we have the usual H\"older inequality 
$$\int_{K}|f(x)g(x)|dx\leq \left(\int_{K}|f(x)|^pdx\right)^{\frac{1}{p}} \left(\int_{K}|f(x)|^{p'}dx\right)^{\frac{1}{p'}},$$
for any functions $f\in L^p(K)$ and $g\in L^{p'}(K)$.\\

See more details about function spaces in the framework of $\mathfrak{p}$-adic spaces in \cite{Kim1}, \cite{Kim2} and \cite{VVZ} as well as the references therein.\\

\mysection{Proof of the Theorem \ref{Theorem1}}\label{ProofTheorem1}
Associated to the operator $T_\Omega$ given in the expression (\ref{Def_Operator_intro}) we can consider the maximal operator $T^*_{\Omega}$ given by 
$$T^*_{\Omega}(f)(x)=\underset{\epsilon>0}{\sup}\left|\int_{\{\epsilon<\|x-y\|_{\mathfrak{p}}\}}\Omega(x,y)f(y)dy\right|,$$
and for some $\epsilon>0$ we consider now the truncated operator
$$T^\epsilon_{\Omega}(f)(x)=\int_{\{\epsilon<\|x-y\|_{\mathfrak{p}}\}}\Omega(x,y)f(y)dy.$$
We recall that we have $T^*_{\Omega}(f)(x)=\underset{\epsilon>0}{\sup}|T^\epsilon_{\Omega}(f)(x)|$ and that 
$$|T_\Omega(f)(x)|\leq T^*_{\Omega}(f)(x),$$
let us remark that the previous control will be very useful in the sequel.\\

Now, for a regular function $f\in L^1_{loc}(\mathbb{Q}_{\mathfrak{p}})$ and for some $k_0\in \mathbb{Z}$ so that $\mathfrak{p}^{k_0-2}<\epsilon < \mathfrak{p}^{k_0-1}$, we write
$$T^\epsilon_{\Omega}(f)(x)=\int_{\{\epsilon<\|x-y\|_{\mathfrak{p}}\leq \mathfrak{p}^{k_0-1}\}}\Omega(x,y)f(y)dy+\sum_{k\geq k_0}\int_{\{\mathfrak{p}^{k-1}<\|x-y\|_{\mathfrak{p}}\leq \mathfrak{p}^{k}\}}\Omega(x,y)f(y)dy.$$
At this point we use the cancellation property (\ref{PseudoRadial}) of the kernel $\Omega$ to introduce some constants: 
$$T^\epsilon_{\Omega}(f)(x)=\int_{\{\epsilon<\|x-y\|_{\mathfrak{p}}\leq \mathfrak{p}^{k_0-1}\}}\Omega(x,y)(f(y)-c_k)dy+\sum_{k\geq k_0}\int_{\{\mathfrak{p}^{k-1}<\|x-y\|_{\mathfrak{p}}\leq \mathfrak{p}^{k}\}}\Omega(x,y)(f(y)-c_k)dy,$$
from which we deduce the inequality
$$|T^\epsilon_{\Omega}(f)(x)|\leq\sum_{k\in\mathbb{Z}}\int_{\{\mathfrak{p}^{k-1}<\|x-y\|_{\mathfrak{p}}\leq \mathfrak{p}^{k}\}}|\Omega(x,y)|  |f(y)-c_k|dy.$$
Now, since by the hypothesis (\ref{SizeKernel}) we have the control $|\Omega(x,y)|\leq \frac{C_\Omega}{\|x-y\|_{\mathfrak{p}}^n}$ for all $x\neq y$, we can write (recalling that we integrate over the set $\{\mathfrak{p}^{k-1}<\|x-y\|_{\mathfrak{p}}^n\leq \mathfrak{p}^{k}\}$)
\begin{eqnarray*}
|T^\epsilon_{\Omega}(f)(x)|&\leq& \sum_{k\in\mathbb{Z}}\int_{\{\mathfrak{p}^{k-1}<\|x-y\|_{\mathfrak{p}}^n\leq \mathfrak{p}^{k}\}}\frac{C_\Omega}{\|x-y\|_{\mathfrak{p}}^n}|f(y)-c_k|dy\\
&\leq &C_\Omega\sum_{k\in\mathbb{Z}}\int_{\{\mathfrak{p}^{k-1}<\|x-y\|_{\mathfrak{p}}^n\leq \mathfrak{p}^{k}\}}\frac{1}{(\mathfrak{p}^{k-1})^n}|f(y)-c_k|dy,
\end{eqnarray*}
 which we rewrite as follows (since $\mathfrak{p}^{nk}=|B(x, \mathfrak{p}^k)|$)
 \begin{eqnarray*}
|T^\epsilon_{\Omega}(f)(x)|&\leq &C_\Omega\cdot \mathfrak{p}^n\sum_{k\in\mathbb{Z}}\int_{\{\|x-y\|_{\mathfrak{p}}^n\leq \mathfrak{p}^{k}\}}\frac{1}{\mathfrak{p}^{nk}}|f(y)-c_k|dy\\
&\leq &C_\Omega\cdot \mathfrak{p}^n\sum_{k\in\mathbb{Z}}\frac{1}{|B(x, \mathfrak{p}^k)|}\int_{\{\|x-y\|_{\mathfrak{p}}^n\leq \mathfrak{p}^{k}\}}|f(y)-c_k|dy\\
&\leq & C_\Omega\cdot \mathfrak{p}^n\sum_{k\in\mathbb{Z}}\frac{1}{|B(x, \mathfrak{p}^k)|}\int_{B(x, \mathfrak{p}^k)}|f(y)-c_k|dy.
\end{eqnarray*}
We fix now the constants $c_k$ by setting $c_k=\frac{1}{|B(x, \mathfrak{p}^k)|}\displaystyle{\int_{B(x, \mathfrak{p}^k)}f(y)dy}=f_B$ and we obtain the inequality
$$|T^\epsilon_{\Omega}(f)(x)|\leq C_\Omega\cdot \mathfrak{p}^n\sum_{k\in\mathbb{Z}}\frac{1}{|B(x, \mathfrak{p}^k)|}\int_{B(x, \mathfrak{p}^k)}|f(y)-f_B|dy.$$
At this point, we decompose the previous sum in two parts by introducing a threshold $\mathcal{K}>0$ to be fixed later:  
\begin{equation}\label{Somme1}
\begin{split}
|T^\epsilon_{\Omega}(f)(x)|&\leq C \underbrace{\sum_{k\leq \lfloor \log_\mathfrak{p}(\mathcal{K})\rfloor}\frac{1}{|B(x, \mathfrak{p}^k)|}\int_{B(x, \mathfrak{p}^k)}|f(y)-f_B|dy}_{(i)}\\
&+ C\underbrace{\sum_{k>\lfloor \log_\mathfrak{p}(\mathcal{K})\rfloor}\frac{1}{|B(x, \mathfrak{p}^k)|}\int_{B(x, \mathfrak{p}^k)}|f(y)-f_B|dy}_{(ii)},
\end{split}
\end{equation}
and we will estimate each of the previous terms separately. 
\begin{itemize}
\item For the first term above, we have, for some index $0<\sigma<1$:
\begin{eqnarray*}
(i)&=&\sum_{k\leq \lfloor \log_\mathfrak{p}(\mathcal{K})\rfloor}\frac{1}{|B(x, \mathfrak{p}^k)|}\int_{B(x, \mathfrak{p}^k)}|f(y)-f_B|dy\\
&=&\sum_{k\leq \lfloor \log_\mathfrak{p}(\mathcal{K})\rfloor}\frac{\mathfrak{p}^{\sigma k}}{(\mathfrak{p}^{k})^{(n+\sigma)}}\int_{B(x, \mathfrak{p}^k)}|f(y)-f_B|dy,
\end{eqnarray*}
and by the definition of the Morrey-Campanato spaces given in the expression (\ref{Def_MorreyCampanato}) above, we obtain:
$$(i)\leq \underset{r>0}{\sup}\underset{x\in \mathbb{Q}^n_{\mathfrak{p}}}{\sup}\left(\frac{1}{r^{(n+\sigma)}}\int_{B(x, r)}|f(y)-f_B|dy\right)\sum_{k\leq \lfloor \log_\mathfrak{p}(\mathcal{K})\rfloor}\mathfrak{p}^{\sigma k}.$$
Since the last sum converges, we obtain 
$$(i)\leq  C\|f\|_{\mathcal{M}^{1, n+\sigma}}\mathcal{K}^{\sigma}.$$
\item For the term $(ii)$ in (\ref{Somme1}) we consider the following quantity for some parameter $0<\alpha<n$ 
\begin{eqnarray*}
(ii)&=&\sum_{k>\lfloor \log_\mathfrak{p}(\mathcal{K})\rfloor}\frac{1}{|B(x, \mathfrak{p}^k)|}\int_{B(x, \mathfrak{p}^k)}|f(y)-f_B|dy\\
&=&\sum_{k>\lfloor \log_\mathfrak{p}(\mathcal{K})\rfloor}\mathfrak{p}^{-k (n\frac{\alpha}{n})}\frac{\mathfrak{p}^{k (n\frac{\alpha}{n})}}{|B(x, \mathfrak{p}^k)|}\int_{B(x, \mathfrak{p}^k)}|f(y)-f_B|dy.
\end{eqnarray*}
Since we have $f_B=\frac{1}{|B(x, \mathfrak{p}^k)|}\displaystyle{\int_{B(x, \mathfrak{p}^k)}f(y)dy}$, we have
\begin{eqnarray*}
(ii)&= &\sum_{k>\lfloor \log_\mathfrak{p}(\mathcal{K})\rfloor}\mathfrak{p}^{-k \alpha}\frac{\mathfrak{p}^{k \alpha}}{|B(x, \mathfrak{p}^k)|}\left[\int_{B(x, \mathfrak{p}^k)}\left|f(y)-\left(\frac{1}{|B(x, \mathfrak{p}^k)|}\displaystyle{\int_{B(x, \mathfrak{p}^k)}f(y)dy}\right)\right|dy\right]\\
&\leq & \sum_{k>\lfloor \log_\mathfrak{p}(\mathcal{K})\rfloor}\mathfrak{p}^{-k \alpha}\frac{\mathfrak{p}^{k \alpha}}{|B(x, \mathfrak{p}^k)|}\left[\int_{B(x, \mathfrak{p}^k)}|f(y)|dy+\int_{B(x, \mathfrak{p}^k)}|f(y)|dy\right]\\
&\leq & 2\sum_{k>\lfloor \log_\mathfrak{p}(\mathcal{K})\rfloor}\mathfrak{p}^{-k \alpha}\frac{\mathfrak{p}^{k \alpha}}{|B(x, \mathfrak{p}^k)|}\int_{B(x, \mathfrak{p}^k)}|f(y)|dy.
\end{eqnarray*}
Recalling now the definition of the fractional maximal function $\mathscr{M}_{\alpha}$ given in the formula (\ref{Def_FuncMax}), we can write 
\begin{eqnarray*}
(ii)&\leq &C \underset{\mathfrak{p}^k>0}{\sup}\left(\frac{\mathfrak{p}^{k \alpha}}{|B(x, \mathfrak{p}^k)|}\int_{B(x, \mathfrak{p}^k)}|f(y)|dy\right)\sum_{k>\lfloor \log_\mathfrak{p}(\mathcal{K})\rfloor}\mathfrak{p}^{-k \alpha}\\
&\leq & C\mathscr{M}_{\alpha}(f)(x)\sum_{k>\lfloor \log_\mathfrak{p}(\mathcal{K})\rfloor}\mathfrak{p}^{-k \alpha}.
\end{eqnarray*}
Noting that the previous sum converges, we finally obtain the control 
$$(ii)\leq C\mathscr{M}_{\alpha}(f)(x) \mathcal{K}^{-\alpha}.\\$$
\end{itemize}
\noindent With the controls obtained for $(i)$ and $(ii)$ we can come back to the expression (\ref{Somme1}) and we have 
$$|T^\epsilon_{\Omega}(f)(x)|\leq (i)+(ii)\leq C\|f\|_{\mathcal{M}^{1, n+\sigma}}\mathcal{K}^{\sigma}+C\mathscr{M}_{\alpha}(f)(x) \mathcal{K}^{-\alpha}.$$
If we set now the constant $\mathcal{K}$ by writing
$$\mathcal{K}=\left(\frac{\mathscr{M}_{\alpha}(f)}{\|f\|_{\mathcal{M}^{1, n+\sigma}}}\right)^{\frac{1}{\alpha+\sigma}},$$
we obtain the estimate
$$|T^\epsilon_{\Omega}(f)(x)|\leq C\mathscr{M}_{\alpha}(f)(x)^{\frac{\sigma}{\alpha+\sigma}}\|f\|_{\mathcal{M}^{1, n+\sigma}}^{1-\frac{\sigma}{\alpha+\sigma}},$$
and from this inequality (which is uniform in $\epsilon>0$) we deduce the wished control
$$T^*_{\Omega}(f)(x)\leq C\mathscr{M}_{\alpha}(f)(x)^{\frac{\sigma}{\alpha+\sigma}}\|f\|_{\mathcal{M}^{1, n+\sigma}}^{1-\frac{\sigma}{\alpha+\sigma}}.$$
Now, setting $\theta=\frac{\sigma}{\alpha+\sigma}$, we have
$$T^*_{\Omega}(f)(x)\leq C\mathscr{M}_{\alpha}(f)(x)^{\theta}\|f\|_{\mathcal{M}^{1, n+\sigma}}^{1-\theta},$$
from which we finally deduce the wished estimate
$$|T_{\Omega}(f)(x)|\leq C\mathscr{M}_{\alpha}(f)(x)^{\theta}\|f\|_{\mathcal{M}^{1, n+\sigma}}^{1-\theta},$$
and this finishes the proof of the Theorem \ref{Theorem1}. \hfill $\blacksquare$\\
\mysection{Proof of the Theorem \ref{Theorem2}}\label{ProofTheorem2}
We note that the proof of the Theorem \ref{Theorem2} follows, in its main lines, the one of the Theorem \ref{Theorem1} presented above. We will give here the details for the sake of completeness and we start by considering the truncated operator
$$T^\epsilon_{\Omega}(f)(x)=\int_{\{\epsilon<\|x-y\|_{\mathfrak{p}}\}}\frac{\omega(x-y)}{\|x-y\|_{\mathfrak{p}}^n}f(y)dy,$$
for some $0<\epsilon<+\infty$ and for some $k_0\in \mathbb{Z}$ so that $\mathfrak{p}^{k_0-2}<\epsilon < \mathfrak{p}^{k_0-1}$, we write
$$T^\epsilon_{\Omega}(f)(x)=\int_{\{\epsilon<\|x-y\|_{\mathfrak{p}}\leq \mathfrak{p}^{k_0-1}\}}\frac{\omega(x-y)}{\|x-y\|_{\mathfrak{p}}^n}f(y)dy+\sum_{k\geq k_0}\int_{\{\mathfrak{p}^{k-1}<\|x-y\|_{\mathfrak{p}}\leq \mathfrak{p}^{k}\}}\frac{\omega(x-y)}{\|x-y\|_{\mathfrak{p}}^n}f(y)dy.$$
At this point we use the cancellation property (\ref{PseudoRadial1}) of the function $\frac{\omega(z)}{\|z\|_{\mathfrak{p}}^n}$ to introduce some constants: 
\begin{eqnarray*}
T^\epsilon_{\Omega}(f)(x)&=&\int_{\{\epsilon<\|x-y\|_{\mathfrak{p}}\leq \mathfrak{p}^{k_0-1}\}}\frac{\omega(x-y)}{\|x-y\|_{\mathfrak{p}}^n}(f(y)-c_k)dy\\
&&+\sum_{k\geq k_0}\int_{\{\mathfrak{p}^{k-1}<\|x-y\|_{\mathfrak{p}}\leq \mathfrak{p}^{k}\}}\frac{\omega(x-y)}{\|x-y\|_{\mathfrak{p}}^n}(f(y)-c_k)dy,
\end{eqnarray*}
and we obtain the inequality
$$|T^\epsilon_{\Omega}(f)(x)|\leq\sum_{k\in\mathbb{Z}}\int_{\{\mathfrak{p}^{k-1}<\|x-y\|_{\mathfrak{p}}\leq \mathfrak{p}^{k}\}}\frac{|\omega(x-y)|}{\|x-y\|_{\mathfrak{p}}^n} |f(y)-c_k|dy.$$
We set now $c_k=\frac{1}{|B(x, \mathfrak{p}^k)|}\displaystyle{\int_{B(x, \mathfrak{p}^k)}f(y)dy}=f_B$ and we introduce a threshold $\mathscr{K}>0$ to be fixed later:  
\begin{equation}\label{Somme2}
\begin{split}
|T^\epsilon_{\Omega}(f)(x)|&\leq \underbrace{\sum_{k\leq \lfloor \log_\mathfrak{p}(\mathscr{K})\rfloor}\int_{\{\mathfrak{p}^{k-1}<\|x-y\|_{\mathfrak{p}}\leq \mathfrak{p}^{k}\}}\frac{|\omega(x-y)|}{\|x-y\|_{\mathfrak{p}}^n} |f(y)-f_B|dy}_{(a)}\\
&+\underbrace{\sum_{k>\lfloor \log_\mathfrak{p}(\mathscr{K})\rfloor}\int_{\{\mathfrak{p}^{k-1}<\|x-y\|_{\mathfrak{p}}\leq \mathfrak{p}^{k}\}}\frac{|\omega(x-y)|}{\|x-y\|_{\mathfrak{p}}^n} |f(y)-f_B|dy}_{(b)},
\end{split}
\end{equation}
and, as before, we will estimate each of the previous terms separately. 
\begin{itemize}
\item For the first term above, we have
\begin{eqnarray*}
(a)&=&\sum_{k\leq \lfloor \log_\mathfrak{p}(\mathscr{K})\rfloor}\int_{\{\mathfrak{p}^{k-1}<\|x-y\|_{\mathfrak{p}}\leq \mathfrak{p}^{k}\}}\frac{|\omega(x-y)|}{\|x-y\|_{\mathfrak{p}}^n} |f(y)-f_B|dy\\
&\leq & C\sum_{k\leq \lfloor \log_\mathfrak{p}(\mathcal{K})\rfloor}\int_{\{\mathfrak{p}^{k-1}<\|x-y\|_{\mathfrak{p}}\leq \mathfrak{p}^{k}\}}\frac{|\omega(x-y)|}{\mathfrak{p}^{kn}} |f(y)-f_B|dy,
\end{eqnarray*}
since over the set $\{\mathfrak{p}^{k-1}<\|x-y\|_{\mathfrak{p}}\leq \mathfrak{p}^{k}\}$ we have $\frac{1}{\|x-y\|_{\mathfrak{p}}^n}\leq \frac{\mathfrak{p}^{n}}{\mathfrak{p}^{kn}}$.
Now, by the H\"older inequality with $1=\frac{1}{\rho}+\frac{1}{\rho'}$, we obtain
\begin{eqnarray*}
(a)\leq C\sum_{k\leq \lfloor \log_\mathfrak{p}(\mathcal{K})\rfloor}\left(\int_{\{\|x-y\|_{\mathfrak{p}}\leq \mathfrak{p}^{k}\}}|\omega(x-y)|^\rho dy\right)^{\frac{1}{\rho}}\left(\int_{\{\|x-y\|_{\mathfrak{p}}\leq \mathfrak{p}^{k}\}}\frac{1}{\mathfrak{p}^{kn\rho'}} |f(y)-f_B|^{\rho'}dy\right)^{\frac{1}{\rho'}}\\
\leq C\sum_{k\leq \lfloor \log_\mathfrak{p}(\mathcal{K})\rfloor}\left(\int_{\mathbb{Q}^n_\mathfrak{p}}|\omega(y)|^\rho dy\right)^{\frac{1}{\rho}}\left(\int_{\{\|x-y\|_{\mathfrak{p}}\leq \mathfrak{p}^{k}\}}\frac{1}{\mathfrak{p}^{kn\rho'}} |f(y)-f_B|^{\rho'}dy\right)^{\frac{1}{\rho'}},
\end{eqnarray*}
which we rewrite as follows since by hypothesis $\|\omega\|_{L^\rho}<+\infty$:
$$(a)\leq C\|\omega\|_{L^\rho}\sum_{k\leq \lfloor \log_\mathfrak{p}(\mathcal{K})\rfloor}\left(\frac{1}{\mathfrak{p}^{kn\rho'}} \int_{B(x,\mathfrak{p}^k)}|f(y)-f_B|^{\rho'}dy\right)^{\frac{1}{\rho'}}.$$
Introducing some power of $\mathfrak{p}$ we obtain
\begin{eqnarray*}
(a)&\leq &C\|\omega\|_{L^\rho}\sum_{k\leq \lfloor \log_\mathfrak{p}(\mathcal{K})\rfloor}\left(\frac{\mathfrak{p}^{k(n+\sigma-n\rho')}}{\mathfrak{p}^{kn\rho'}\times \mathfrak{p}^{k(n+\sigma-n\rho')}} \int_{B(x,\mathfrak{p}^k)}|f(y)-f_B|^{\rho'}dy\right)^{\frac{1}{\rho'}}\\
&\leq &C\|\omega\|_{L^\rho}\sum_{k\leq \lfloor \log_\mathfrak{p}(\mathcal{K})\rfloor}\left(\frac{\mathfrak{p}^{k(n+\sigma-n\rho')}}{\mathfrak{p}^{k(n+\sigma)}} \int_{B(x,\mathfrak{p}^k)}|f(y)-f_B|^{\rho'}dy\right)^{\frac{1}{\rho'}},
\end{eqnarray*}
which gives us the following estimate
$$(a)\leq C\|\omega\|_{L^\rho}\sum_{k\leq \lfloor \log_\mathfrak{p}(\mathcal{K})\rfloor}(\mathfrak{p}^k)^{\frac{(n+\sigma-n\rho')}{\rho'}}\left(\frac{1}{\mathfrak{p}^{k(n+\sigma)}} \int_{B(x,\mathfrak{p}^k)}|f(y)-f_B|^{\rho'}dy\right)^{\frac{1}{\rho'}}.$$
By the definition of the Morrey-Campanato spaces given in the expression (\ref{Def_MorreyCampanato}) above, we obtain:
\begin{eqnarray*}
(a)&\leq &C\underset{r>0}{\sup}\underset{x\in \mathbb{Q}^n_{\mathfrak{p}}}{\sup}\left(\frac{1}{r^{(n+\sigma)}}\int_{B(x, r)}|f(y)-f_B|^{\rho'}dy\right)^{\frac{1}{\rho'}}\sum_{k\leq \lfloor \log_\mathfrak{p}(\mathcal{K})\rfloor}(\mathfrak{p}^k)^{\frac{(n+\sigma-n\rho')}{\rho'}}\\
&\leq & C\|f\|_{\mathcal{M}^{\rho', n+\sigma}}\sum_{k\leq \lfloor \log_\mathfrak{p}(\mathcal{K})\rfloor}(\mathfrak{p}^k)^{\frac{(n+\sigma-n\rho')}{\rho'}}.
\end{eqnarray*}
Recall that by hypothesis we have $1<\rho'<\frac{n+\sigma}{n}$, which guarantees that $n+\sigma-n\rho'>0$, and thus the previous sum converges, se we have 
$$(a)\leq  C\|f\|_{\mathcal{M}^{\rho', n+\sigma}}\mathcal{K}^{\frac{(n+\sigma-n\rho')}{\rho'}}.$$
\item For the term $(b)$ in (\ref{Somme2}) we have
$$(b)=\sum_{k>\lfloor \log_\mathfrak{p}(\mathscr{K})\rfloor}\int_{\{\mathfrak{p}^{k-1}<\|x-y\|_{\mathfrak{p}}\leq \mathfrak{p}^{k}\}}\frac{|\omega(x-y)|}{\|x-y\|_{\mathfrak{p}}^n} |f(y)-f_B|dy,$$
and since over the set $\{\mathfrak{p}^{k-1}<\|x-y\|_{\mathfrak{p}}\leq \mathfrak{p}^{k}\}$ we have $\frac{1}{\|x-y\|_{\mathfrak{p}}^n}\leq \frac{\mathfrak{p}^{n}}{\mathfrak{p}^{kn}}$, we can write
$$(b)\leq C\sum_{k>\lfloor \log_\mathfrak{p}(\mathscr{K})\rfloor}\int_{\{\mathfrak{p}^{k-1}<\|x-y\|_{\mathfrak{p}}\leq \mathfrak{p}^{k}\}}\frac{|\omega(x-y)|}{\mathfrak{p}^{nk}} |f(y)-f_B|dy.$$
Now, by the H\"older inequality with $1=\frac{1}{\rho}+\frac{1}{\rho'}$, we have
\begin{eqnarray*}
(b)&\leq &C\|\omega\|_{L^{\rho}}\sum_{k>\lfloor \log_\mathfrak{p}(\mathscr{K})\rfloor}\frac{1}{\mathfrak{p}^{nk}}\left(\int_{\{\mathfrak{p}^{k-1}<\|x-y\|_{\mathfrak{p}}\leq \mathfrak{p}^{k}\}} |f(y)-f_B|^{\rho'}dy\right)^{\frac{1}{\rho'}}\\
&\leq & C\|\omega\|_{L^{\rho}}\sum_{k>\lfloor \log_\mathfrak{p}(\mathscr{K})\rfloor}\frac{1}{\mathfrak{p}^{nk}}\left(\int_{\{\|x-y\|_{\mathfrak{p}}\leq \mathfrak{p}^{k}\}} |f(y)-f_B|^{\rho'}dy\right)^{\frac{1}{\rho'}},
\end{eqnarray*}
from which we deduce, by the triangular inequality in $L^{\rho'}$ spaces, the control 
\begin{eqnarray*}
(b)&\leq & C\|\omega\|_{L^{\rho}}\sum_{k>\lfloor \log_\mathfrak{p}(\mathscr{K})\rfloor}\frac{1}{\mathfrak{p}^{nk}}\left[\left(\int_{B(x, \mathfrak{p}^k)} |f(y)|^{\rho'}dy\right)^{\frac{1}{\rho'}}+\left(\int_{B(x, \mathfrak{p}^k)} |f_B|^{\rho'}dy\right)^{\frac{1}{\rho'}}\right]\\
&\leq & C\|\omega\|_{L^{\rho}}\sum_{k>\lfloor \log_\mathfrak{p}(\mathscr{K})\rfloor}\frac{1}{\mathfrak{p}^{nk}}\left[\left(\int_{B(x, \mathfrak{p}^k)} |f(y)|^{\rho'}dy\right)^{\frac{1}{\rho'}}+|f_B| |B(x, \mathfrak{p}^k)|^{\frac{1}{\rho'}}\right].
\end{eqnarray*}
Recalling the definition of the quantities $f_B$ and since $|B(x), \mathfrak{p}^k|^{\frac{1}{\rho'}}=\mathfrak{p}^{\frac{nk}{\rho'}}$, we have now
\begin{eqnarray*}
(b)&\leq & C\sum_{k>\lfloor \log_\mathfrak{p}(\mathscr{K})\rfloor}\frac{1}{\mathfrak{p}^{nk}}\left[\left(\int_{B(x, \mathfrak{p}^k)} |f(y)|^{\rho'}dy\right)^{\frac{1}{\rho'}}+\left|\frac{1}{|B(x, \mathfrak{p}^k)|}\int_{B(x, \mathfrak{p}^k)}f(y)dy\right| \mathfrak{p}^{\frac{nk}{\rho'}}\right]\\
&\leq & C\sum_{k>\lfloor \log_\mathfrak{p}(\mathscr{K})\rfloor}\frac{1}{\mathfrak{p}^{nk}}\left[\left(\int_{B(x, \mathfrak{p}^k)} |f(y)|^{\rho'}dy\right)^{\frac{1}{\rho'}}+\frac{|B(x, \mathfrak{p}^k)|^{\frac{1}{\rho}}}{\mathfrak{p}^{nk}}\left(\int_{B(x, \mathfrak{p}^k)} |f(y)|^{\rho'}dy\right)^{\frac{1}{\rho'}} \mathfrak{p}^{\frac{nk}{\rho'}}\right],
\end{eqnarray*}
where we applied the H\"older inequality with $1=\frac{1}{\rho}+\frac{1}{\rho'}$ in the last integral above, and we thus obtain
\begin{eqnarray*}
(b)&\leq & C\sum_{k>\lfloor \log_\mathfrak{p}(\mathscr{K})\rfloor}\frac{1}{\mathfrak{p}^{nk}}\left[\left(\int_{B(x, \mathfrak{p}^k)} |f(y)|^{\rho'}dy\right)^{\frac{1}{\rho'}}+\frac{\mathfrak{p}^{\frac{nk}{\rho}}}{\mathfrak{p}^{nk}}\left(\int_{B(x, \mathfrak{p}^k)} |f(y)|^{\rho'}dy\right)^{\frac{1}{\rho'}} \mathfrak{p}^{\frac{nk}{\rho'}}\right]\\
&\leq & C'\sum_{k>\lfloor \log_\mathfrak{p}(\mathscr{K})\rfloor}\frac{1}{\mathfrak{p}^{nk}}\left(\int_{B(x, \mathfrak{p}^k)} |f(y)|^{\rho'}dy\right)^{\frac{1}{\rho'}}.
\end{eqnarray*}
We write now
\begin{eqnarray*}
(b)&\leq &C\sum_{k>\lfloor \log_\mathfrak{p}(\mathscr{K})\rfloor}\frac{1}{\mathfrak{p}^{nk}}\left(\frac{\mathfrak{p}^{nk}}{\mathfrak{p}^{\alpha k}}\right)^{\frac{1}{\rho'}}\left(\frac{\mathfrak{p}^{\alpha k}}{\mathfrak{p}^{nk}}\int_{B(x, \mathfrak{p}^k)} |f(y)|^{\rho'}dy\right)^{\frac{1}{\rho'}}\\
&\leq &C\sum_{k>\lfloor \log_\mathfrak{p}(\mathscr{K})\rfloor}\mathfrak{p}^{k(\frac{n}{\rho'}-n-\frac{\alpha}{\rho'})}\left(\frac{\mathfrak{p}^{\alpha k}}{|B(x, \mathfrak{p}^k)|}\int_{B(x, \mathfrak{p}^k)} |f(y)|^{\rho'}dy\right)^{\frac{1}{\rho'}}.
\end{eqnarray*}
Recalling now the definition of the fractional maximal function $\mathscr{M}_{\alpha}$ given in the formula (\ref{Def_FuncMax}), we can write 
\begin{eqnarray*}
(b)&\leq &C \left[\underset{\mathfrak{p}^k>0}{\sup}\left(\frac{\mathfrak{p}^{k \alpha}}{|B(x, \mathfrak{p}^k)|}\int_{B(x, \mathfrak{p}^k)}|f(y)|^{\rho'}dy\right)\right]^{\frac{1}{\rho'}}\sum_{k>\lfloor \log_\mathfrak{p}(\mathcal{K})\rfloor}(\mathfrak{p}^{k})^{-\frac{1}{\rho}-\frac{\alpha}{\rho'}}\\
&\leq & C[\mathscr{M}_{\alpha}(|f|^{\rho'})(x)]^{\frac{1}{\rho'}}\sum_{k>\lfloor \log_\mathfrak{p}(\mathcal{K})\rfloor}(\mathfrak{p}^{k})^{-\frac{1}{\rho}-\frac{\alpha}{\rho'}},
\end{eqnarray*}
and since the previous sum converges, we finally obtain the control 
$$(b)\leq C[\mathscr{M}_{\alpha}(|f|^{\rho'})(x)]^{\frac{1}{\rho'}} \mathcal{K}^{-\frac{1}{\rho}-\frac{\alpha}{\rho'}}.\\$$
\end{itemize}
\noindent With the controls obtained for $(a)$ and $(b)$ we can come back to the expression (\ref{Somme2}) and we have 
$$|T^\epsilon_{\Omega}(f)(x)|\leq (a)+(b)\leq C\|f\|_{\mathcal{M}^{\rho', n+\sigma}}\mathcal{K}^{\frac{(n+\sigma-n\rho')}{\rho'}}+C[\mathscr{M}_{\alpha}(|f|^{\rho'})(x)]^{\frac{1}{\rho'}} \mathcal{K}^{-\frac{1}{\rho}-\frac{\alpha}{\rho'}}.$$
If we set now the constant $\mathcal{K}$ by 
$$\mathcal{K}=\left(\frac{[\mathscr{M}_{\alpha}(|f|^{\rho'})]^{\frac{1}{\rho'}}}{\|f\|_{\mathcal{M}^{\rho', n+\sigma}}}\right)^{\frac{1}{\frac{(n+\sigma-n\rho')}{\rho'}+\frac{1}{\rho}+\frac{\alpha}{\rho'}}},$$
we obtain the estimate
$$|T^\epsilon_{\Omega}(f)(x)|\leq C[\mathscr{M}_{\alpha}(|f|^{\rho'})(x)]^{\frac{1}{\rho'}(\frac{\frac{(n+\sigma-n\rho')}{\rho'}}{\frac{(n+\sigma-n\rho')}{\rho'}+\frac{1}{\rho}+\frac{\alpha}{\rho'}})}\|f\|_{\mathcal{M}^{\rho', n+\sigma}}^{\frac{\frac{1}{\rho}+\frac{\alpha}{\rho'}}{\frac{(n+\sigma-n\rho')}{\rho'}+\frac{1}{\rho}+\frac{\alpha}{\rho'}}}.$$
If we set $0<\theta=\frac{\frac{(n+\sigma-n\rho')}{\rho'}}{\frac{(n+\sigma-n\rho')}{\rho'}+\frac{1}{\rho}+\frac{\alpha}{\rho'}}<1$, we have
$$|T^\epsilon_{\Omega}(f)(x)|\leq C[\mathscr{M}_{\alpha}(|f|^{\rho'})(x)]^{\frac{\theta}{\rho'}}\|f\|_{\mathcal{M}^{\rho', n+\sigma}}^{1-\theta}.$$
and from this inequality (which is uniform in $\epsilon>0$) we deduce the wished control
$$T^*_{\Omega}(f)(x)\leq C[\mathscr{M}_{\alpha}(|f|^{\rho'})(x)]^{\frac{\theta}{\rho'}}\|f\|_{\mathcal{M}^{\rho', n+\sigma}}^{1-\theta},$$
\emph{i.e.}: 
$$|T_{\Omega}(f)(x)|\leq C[\mathscr{M}_{\alpha}(|f|^{\rho'})(x)]^{\frac{\theta}{\rho'}}\|f\|_{\mathcal{M}^{\rho', n+\sigma}}^{1-\theta},$$
and this finishes the proof of the Theorem \ref{Theorem2}. \hfill $\blacksquare$\\
\mysection{Proof of the Corollary \ref{Coro1}}\label{ProofCoro1}
In order to proof the Corollary \ref{Coro1}, we first observe that we have the following pointwise control
\begin{equation}\label{MaxFracRiesz}
\mathscr{M}_\alpha(f)(x)\leq I_\alpha(|f|)(x),
\end{equation}
where 
\begin{equation}\label{DefRiesz}
I_{\alpha}(f)(x)=\int_{\mathbb{Q}^n_\mathfrak{p}}\frac{|f(y)|}{\|x-y\|^{n-\alpha}_\mathfrak{p}}dy, 
\end{equation}
is the usual Riesz potential\footnote{Note that in the literature a different normalization is used: a factor $\frac{1}{\Gamma_n(\alpha)}$ with $\Gamma_n(\alpha)=\frac{1-\mathfrak{p}^{-n-\alpha}}{1-\mathfrak{p}^{-\alpha}}$ is usually added, see \emph{e.g.} \cite{Shi} and the references therein for more details.}. Indeed, for some fixed radius $r>0$ we write
$$\frac{r^\alpha}{|B(x,r)|}\int_{B(x,r)}|f(y)|dy=\frac{r^\alpha}{r^n}\int_{\{\|x-y\|_{\mathfrak{p}}\leq r\}}|f(y)|dy=\int_{\{\|x-y\|_{\mathfrak{p}}\leq r\}}\frac{|f(y)|}{r^{n-\alpha}}dy.$$
Since over the set $\{\|x-y\|_{\mathfrak{p}}\leq r\}$ we have $\frac{1}{r^{n-\alpha}}\leq \frac{1}{\|x-y\|^{n-\alpha}_{\mathfrak{p}}}$, we can write 
$$\frac{r^\alpha}{|B(x,r)|}\int_{B(x,r)}|f(y)|dy\leq \int_{\{\|x-y\|_{\mathfrak{p}}\leq r\}}\frac{|f(y)|}{\|x-y\|^{n-\alpha}_{\mathfrak{p}}}dy\leq  \int_{\mathbb{Q}^n_{\mathfrak{p}}}\frac{|f(y)|}{\|x-y\|^{n-\alpha}_{\mathfrak{p}}}dy=I_{\alpha}(f)(x),$$
and from this uniform estimate (with respect to the radius $r>0$), by the definition of the fractional maximal function $\mathscr{M}_\alpha$ given in (\ref{Def_FuncMax}), we easily deduce the announced control (\ref{MaxFracRiesz}) above. \\

Next, we recall that we have the following Hardy-Littlewood-Sobolev type inequalities in the setting of the $\mathfrak{p}$-adic groups:
\begin{equation}\label{HardyLittlewoodIneq}
\|I_{\alpha}(|f|)\|_{L^q}\leq C\|f\|_{L^p}, \quad \mbox{with }\frac{1}{q}=\frac{1}{p}-\frac{\alpha}{n} \quad\mbox{and } 0<\alpha<n.  
\end{equation}
See a proof of this fact in \cite{Phillips} and in \cite{KimC}.\\

With the estimates (\ref{MaxFracRiesz}) and (\ref{HardyLittlewoodIneq}) at hand, we can now prove the Corollary \ref{Coro1}. Our starting point is the estimate (\ref{PointwiseControl1}):
$$|T_{\Omega}(f)(x)|\leq C[\mathscr{M}_{\alpha}(|f|^{\rho'})(x)]^{\frac{\theta}{\rho'}}\|f\|_{\mathcal{M}^{\rho', n+\sigma}}^{1-\theta}.$$
Using the pointwise control of the fractional maximal function $\mathscr{M}_{\alpha}(|f|^{\rho'})(x)$ by the Riesz potential $I_{\alpha}(|f|^{\rho'})(x)$, we can write 
$$|T_{\Omega}(f)(x)|\leq C[I_{\alpha}(|f|^{\rho'})(x)]^{\frac{\theta}{\rho'}}\|f\|_{\mathcal{M}^{\rho', n+\sigma}}^{1-\theta}.$$
Now we apply an $L^q$-norm to the both sides of this inequality to obtain
$$\|T_{\Omega}(f)\|_{L^q}\leq C\left\|[I_{\alpha}(|f|^{\rho'})]^{\frac{\theta}{\rho'}}\right\|_{L^q}\|f\|_{\mathcal{M}^{\rho', n+\sigma}}^{1-\theta}= C\left\|I_{\alpha}(|f|^{\rho'})\right\|_{L^{\frac{\theta}{\rho'}q}}^{\frac{\theta}{\rho'}}\|f\|_{\mathcal{M}^{\rho', n+\sigma}}^{1-\theta},$$
where we used the identity $\||\varphi|^\beta\|_{L^q}=\|\varphi\|^{\beta}_{L^{\beta q}}$. We use now the Hardy-Littlewood-Sobolev inequality (\ref{HardyLittlewoodIneq}) and we have the estimate $\left\|I_{\alpha}(|f|^{\rho'})\right\|_{L^{\frac{\theta}{\rho'}q}}\leq C\left\||f|^{\rho'}\right\|_{L^{\frac{nq\theta}{n\rho'+\alpha q \theta}}}$ (note that $\frac{\rho'}{q\theta}=\frac{n\rho'+\alpha q \theta}{nq\theta}-\frac{\alpha}{n}$), which allows us to write
$$\|T_{\Omega}(f)\|_{L^q}\leq C\left\||f|^{\rho'}\right\|_{L^{\frac{nq\theta}{n\rho'+\alpha q \theta}}}^{\frac{\theta}{\rho'}}\|f\|_{\mathcal{M}^{\rho', n+\sigma}}^{1-\theta}=C\|f\|_{L^{\frac{n\rho'q\theta}{n\rho'+\alpha q \theta}}}^{\theta}\|f\|_{\mathcal{M}^{\rho', n+\sigma}}^{1-\theta},$$
and since $p=\frac{n\rho'q\theta}{n\rho'+\alpha q \theta}>1$ by hypothesis, we finally obtain 
$$\|T_{\Omega}(f)\|_{L^q}\leq C\|f\|_{L^{p}}^{\theta}\|f\|_{\mathcal{M}^{\rho', n+\sigma}}^{1-\theta},$$
and the proof of the Corollary \ref{Coro1} is now finished. \hfill $\blacksquare$\\

To finish this article, and inspired by the pointwise control (\ref{MaxFracRiesz}), we present now a different version of the Theorem \ref{Theorem2}, that can be useful in some applications, where we replace the fractional maximal function $\mathscr{M}_\alpha$ by the Riesz potential $I_\alpha$:\\

\begin{Corollaire}\label{Coro2}
Under the same hypotheses of the Theorem \ref{Theorem2} for the operator $T_\Omega$, let us  consider $f:\mathbb{Q}^n_{\mathfrak{p}}\longrightarrow \mathbb{R}$ a locally integrable function such that $f$ belongs to the Morrey-Campanato space $\mathcal{M}^{\rho',n+\sigma}(\mathbb{Q}^n_{\mathfrak{p}})$, with $1\leq \rho'\leq n+\sigma<+\infty$ and with $0<\sigma< \rho'$.\\

\noindent If we assume moreover the constraint $\omega\in L^{\rho}(\mathbb{Q}^n_\mathfrak{p})$ with $\frac{n+\sigma}{n+\sigma-1}\leq \rho\leq +\infty$, then we have the pointwise bound 
$$|T_{\Omega}(f)(x)|\leq C[I_{\alpha}(|f|^{\rho'})(x)]^{\frac{\theta}{\rho'}}\|f\|_{\mathcal{M}^{\rho', n+\sigma}}^{1-\theta},$$
where $I_\alpha$ with $0<\alpha<n$ is the Riesz potential given in (\ref{DefRiesz}) and where $0<\theta<1$ is such that $\theta=\frac{\frac{(n+\sigma-n\rho')}{\rho'}}{\frac{(n+\sigma-n\rho')}{\rho'}+\frac{1}{\rho}+\frac{\alpha}{\rho'}}$.\\
\end{Corollaire}
\noindent The proof of this corollary is exactly the same as the one presented above for the Theorem \ref{Theorem2}: it is enough to use the pointwise control (\ref{MaxFracRiesz}) to replace the fractional maximal function $\mathscr{M}_\alpha$ by the Riesz potential $I_\alpha$.\\

\noindent {\bf Conflict of interest.} The author declares no conflicts of interest. \\
\noindent {\bf Data availability statement.} Data sharing is not applicable to this article as no new data were created or analyzed in this study.\\

\end{document}